\documentclass{amsart}
\usepackage{booktabs} 
\usepackage{diagbox}
\newtheorem{theorem}{Theorem}[section]

\theoremstyle{definition}

\newtheorem{example}[theorem]{Example}

\newtheorem{lm}[theorem]{Lemma}

\theoremstyle{remark}

\numberwithin{equation}{section}

\begin{document}

\title{On Characters of a Class of P-polynomial Table Algebras and  Applications}

\author{Masoumeh Koohestani}
\address{Faculty of Mathematics, K. N. Toosi
University of Technology, P. O. Box $16315$-$1618$, Tehran, Iran}

\email{m.kuhestani@email.kntu.ac.ir}
\author{Amir Rahnamai Barghi}
\address{Faculty of Mathematics, K. N. Toosi
University of Technology, P. O. Box $16315$-$1618$, Tehran, Iran}
\email{rahnama@kntu.ac.ir}
\author{Amirhossein Amiraslani}
\address{STEM Department, University of Hawaii-Maui College, Kahului, HI 96732, USA}
\email{aamirasl@hawaii.edu}







\begin{abstract}
In this paper, we study the characters of homogeneous monotonic P-polynomial table algebras with finite dimension $d\geq 5$. We then apply them to association schemes. To this end, we develop some methods using tridiagonal matrices and $\mathbf{z}$-transform. Moreover, we calculate the eigenvalues of a special tridiagonal matrix which is found through the first intersection matrix of P-polynomial table algebras.
\end{abstract}
\maketitle
\section{Introduction}

The characters of table algebras are applied to study the properties of table algebras and  can be used in association schemes and finite groups, see \cite{structure} and \cite{characters}. In particular, calculating the characters of table algebras can help to determine the algebraic structure of association schemes because the eigenvalues of association schemes which determine the algebraic structure of association schemes will be obtained using the characters of table algebras, see \cite{Gods}.
The characters of certain table algebras have also been calculated  in some articles such as \cite{class} for lower dimensions. However, calculating the characters of table algebras is  hard or impossible, in general.

In \cite{table}, the first intersection matrix  and the multiplicities of homogeneous monotonic P-polynomial table algebras are studied, but the characters of this class of table algebras have not been explicitly calculated  yet. Moreover,
the first intersection matrix of P-polynomial table algebras is  tridiagonal and  by calculating the eigenvalues of  the first intersection matrix, we are able to calculate all  characters, see {\cite[Remark $3.1$]{table}}.
 Here, we calculate the characters of homogeneous monotonic P-polynomial table algebras with finite dimension $d\geq 5$ which can be applied in P-polynomial association schemes. To do so, we must study the eigenvalues of a certain tridiagonal matrix.
The eigenstructure and applications  of  tridiagonal matrices have been studied in many articles such as \cite{several}, \cite{system} and
 \cite{Anal}. In this work, we consider a special class
 of tridiagonal matrices that is found through the first intersection  matrix of our desirable table algebras and for which the eigenvalues have not been calculated yet.

The structure of this paper is as follows. In Section~\ref{tablealgebra}, we introduce table algebras and specifically P-polynomial table algebras. Section~\ref{transform1} gives an overview of $\textbf{z}$-transform and some of its properties.  In section \ref{tridiagonal}, we calculate the characteristic polynomial of a tridiagonal matrix which is then used in the following section. In Section~\ref{p-tablealgebra}, we introduce the first intersection matrix of homogeneous monotonic P-polynomial table algebras, and obtain some results regarding its eigenvalues. We also study the characters of the table algebras that we consider in this work. Moreover, we apply our results to study the eigenvalues of two classes of P-polynomial association schemes. Finally, the results of this paper are summarized in Section~\ref{remarks}.

Throughout this paper, $\mathbb{C}$ and $\mathbb{R}^+$ denote the complex numbers and the positive real numbers, respectively.

\section{Table Algebras}\label{tablealgebra}
In this section, we go over some  concepts related to table algebras and P-polynomial table algebras, see~\cite{table} and \cite{scheme}  for more details.

Let $A$ be a finite-dimensional associative commutative algebra with a basis  $\mathbf{B}=\{x_0=1_A, x_1, \cdots , x_d\}$. Then $(A, \mathbf{B})$ is called a table  algebra if the following conditions hold:
\begin{itemize}
\item[(i)]
 $x_ix_j=\sum_{m=0}^d\beta _{ijm}x_m$ with $\beta _{ijm} \in \mathbb{R}^+\cup \{0\}$, for all $i$, $j$;
 \item[(ii)]
 there is an algebra automorphism of $A$ (denoted by $^-$),
whose order divides 2, such that if ${x}_i\in \mathbf{B}$, then $\overline{x}_i\in \mathbf{B}$ and $\overline{i}$ is defined by $x_{\overline{i}}=\overline{x}_i$;
\item[(iii)]
 for all $i$, $j$, we have $\beta_{ij0}\neq 0$ if and only if $j=\overline{i}$, and $\beta_{i\overline{i}0}>0$.
\end{itemize}
Let $(A,\mathbf{B})$ with $\mathbf{B}=\{x_0=1_A, x_1, \cdots , x_d\}$ be  a table algebra. Then $(A,\mathbf{B})$
 is called real, if $i=\overline{i}$, for $0 \leq i\leq d$.
Moreover, the $i$-th intersection matrix of $(A,\mathbf{B})$  (the intersection matrix with respect to $x_i$) is defined as
follows.
$$
B_i=\left(
\begin{array}{cccc}
\beta _{i00} & \beta _{i01} & \ldots & \beta _{i0d}\\
\beta _{i10} & \beta _{i11} & \ldots & \beta _{i1d}\\
\vdots & \vdots & \ddots & \vdots \\
\beta _{id0} & \beta _{id1} & \ldots & \beta _{idd}\\
\end{array} \right)_{(d+1) \times (d+1)},
$$
where $x_ix_j=\sum _{i=0}^d\beta _{ijk}x_k$, for all $i, j, k$.

For any table algebra $(A,\mathbf{B})$ with $\mathbf{B}=\{x_0=1_A, x_1, \cdots , x_d\}$, there exists a unique algebra homomorphism $f:A\rightarrow \mathbb{C}$ such that $f(x_i)=f(x_{\overline{i}})\in \mathbb{R}^+$, for  $0 \leq i\leq d$. If $f(x_i)=\beta_{i\overline{i}0}$, for all $i$, then $(A,\mathbf{B})$ is called standard. If $d\geq 2$ and for $i>0$, $f(x_i)$ is constant and  $(A,\mathbf{B})$ is called homogeneous.

  A  real standard table algebra $(A,\mathbf{B})$ with $\mathbf{B}=\{x_0=1_A, x_1, \cdots , x_d\}$ is called P-polynomial if for each $i$, $2\leq i\leq d$, there exists a complex cofficient polynomial $\nu_i(x)$ of degree $i$ such that $x_i=\nu_i(x_1)$. If $(A,\mathbf{B})$ is a P-polynomial table algebra, then  for all $i$, there exist $b_{i-1}, a_i, c_{i+1}\in \mathbb{R}$ such that
\begin{equation}\label{abc}
x_1x_i=b_{i-1}x_{i-1}+a_ix_i+c_{i+1}x_{i+1},
\end{equation}
 with $b_{i}\neq 0$, ($0\leq i\leq d-1$), $c_{i}\neq 0$, ($1\leq i\leq d$), and  $b_{-1}=c_{d+1}=0$. The first intersection matrix of a P-polynomial table algebra is a tridiagonal matrix as follows.
\begin{equation}\label{matrix}
B_1=\left(
\begin{array}{ccccc}
a_0&c_1&&&\\
b_0&a_1&c_2&&\\
&b_1&a_2&\ddots&\\
&&\ddots &\ddots &c_d\\
&&&b_{d-1}&a_d
\end{array} \right)_{(d+1) \times (d+1)},
\end{equation}
where $b_0=\beta _{110}$ is called the valency of the P-polynomial table algebra. A P-polynomial table algebra is called monotonic if $c_i\leq c_{i+1}$ ($1\leq i\leq d-1$), and $b_i\geq b_{i+1}$ ($0\leq i\leq d-2$).

Let $(A,\mathbf{B})$ with $\mathbf{B}=\{x_0=1_A, x_1, \cdots , x_d\}$ be a table algebra. Since $A$ is semisimple, the primitive idempotents of $A$ form another basis for $A$, see \cite{scheme}.
Consequently, if $\{e_0, e_1, \cdots, e_d\}$ is the set of the primitive idempotents of $A$, then we have $x_i=\sum_{j=0}^dp_i(j)e_j$, where $p_i(j)\in \mathbb{C}$, for $0\leq i\leq d$. The numbers $p_i(j)$ are the characters of the table algebra.

Let $(A,\mathbf{B})$ with $\mathbf{B}=\{x_0=1_A, x_1, \cdots , x_d\}$ be a P-polynomial table algebra. Then the $p_1(j)$ are equal to the eigenvalues of its first intersection matrix and for $2\leq i\leq d$, we have
\begin{equation}\label{character}
p_i(j)=\nu _i(p_1(j)), ~~~ 0\leq j\leq d,
\end{equation}
where $\nu _i(x_1)$ is a complex cofficient polynomial such that $x_i=\nu _i(x_1).$
\section{$\mathbf{z}$-transform}\label{transform1}
This section gives an overview of $\textbf{z}$-transform. The concept of $\textbf{z}$-transform has the same role in discrete-time signals that Laplace transform has in continuous-time signals. For a discrete-time signal which is a sequence of real or complex numbers such as $x[n]$, the $\textbf{z}$-transform is defined as the power series
\begin{equation}\label{transform}
X(z)=\sum _{n=-\infty}^{+\infty}x[n]z^{-n},
\end{equation}
where  $n$ is an integer and $z$ is a complex variable. The function $X(z)$ in (\ref{transform}) is called the two-sided or bilateral $\textbf{z}$-transform of $x[n]$. The one-sided or unilateral $\textbf{z}$-transform of $x[n]$ is defined by
$$X(z)=\sum _{n=0}^{+\infty}x[n]z^{-n}=x[0]+x[1]z^{-1}+x[2]z^{-2}+\cdots.$$
We use the notation $x[n]\leftrightarrow X(z)$ to show that $X(z)$ is the $\textbf{z}$-transform of $x[n]$.  A well-known example of $\textbf{z}$-transform is as follows.
\begin{example}
Let $x[n]$ be a discrete-time function in the form of
$$
x[n]=
\begin{cases}
1,~~n=0, 1, 2, \cdots \\
0,~~n=-1, -2, \cdots
\end{cases}
$$
The $\textbf{z}$-transform $X(z)$ is
\begin{align}\label{x[n]}
X(z)=\sum _{n=0}^{+\infty}x[n]z^{-n}=&\sum _{n=0}^{+\infty}z^{-n}\nonumber\\
=&~1+z^{-1}+z^{-2}+\cdots\nonumber\\
=&~\frac{1}{1-z^{-1}}\nonumber\\
=&~\frac{z}{z-1}.
\end{align}
The region of convergence for the z-transform $X(z)$  in (\ref{x[n]}) is the set of all complex numbers $z$ such that $|z|>1$.
\end{example}
 $\textbf{z}$-transform is a linear operation. This means that if we have $x[n]\leftrightarrow X(z)$ and $u[n]\leftrightarrow U(z)$, then
 \begin{equation}\label{pro1}
 ax[n]+bu[n]\leftrightarrow aX(z)+bU(z),
 \end{equation}
 where  $a, b\in \mathbb{C}$.

 Let $x[n]\leftrightarrow X(z)$ and $q$ be a positive integer. Then  we have
 \begin{equation}\label{pro2}
 x[n-q]\leftrightarrow z^{-q}X(z),
 \end{equation}
 and
 \begin{equation}\label{pro3}
 x[n+q]\leftrightarrow z^qX(z)-x[0]z^q-x[1]z^{q-1}-\cdots -x[q-1]z.
 \end{equation}
The proof of the above properties
and more facts about $\textbf{z}$-transform can be found in
{\cite[Chapter $7$]{matlab}}.

\section{Tridiagonal Matrices}\label{tridiagonal}
    In this section, we find the characteristic polynomial of a tridiagonal matrix, namely $A_n$. The eigenvalues of $A_n$ can be calculated by means of some previous results, e.g. in \cite{eigensystem}. Here, we calculate the characteristic polynomial of $A_n$ through an approach which we apply to study the characters of P-polynomial table algebras in the next section. First, we state the following lemmas:
  \begin{lm}\label{lem1}
(\cite{Fib}) Let $\{H_n, n = 1,2,\cdots \}$ be a sequence of tridiagonal matrices in the form of
$$
H_n= \left(
\begin{array}{ccccc}
h_{1,1}&h_{1,2}&&&\\
h_{2,1}&h_{2,2}&h_{2,3}&0&\\
&h_{3,2}&h_{3,3}&\ddots &\\
&0&\ddots & \ddots &h_{n-1,n}\\
&& &h_{n,n-1} &h_{n,n}\\
\end{array} \right).
$$
Then the determinants of $H_n$ are given by the recursive formula:
\begin{align}
|H_1|&= h_{1,1},\nonumber \\
|H_2|&=h_{1,1}h_{2,2}-h_{1,2}h_{2,1},\nonumber\\
|H_n|&=h_{n,n}|H_{n-1}|-h_{n-1,n}h_{n,n-1}|H_{n-2}|.\nonumber
\end{align}
\end{lm}
\begin{lm}\label{lem2}
Let $A_n$ be the following tridiagonal matrix:
$$
A_n= \left(
\begin{array}{ccccc}
a&1&&&\\
b&a&1&&\\
&b&a&\ddots &\\
&&\ddots & \ddots &1\\
&& &b &a\\
\end{array} \right)_{n \times n},
$$
where $a, b\in \mathbb{C}$ and $b\neq 0$.
Then the characteristic polynomial of $A_n$ is
$$(\sqrt{b})^nU_n\left(\frac{x-a}{2\sqrt{b}}\right),$$
where  $U_n(x)$ is the $n$-th degree Chebyshev polynomial of second kind.

\begin{proof}
Let $\Delta _n(x)$ be a function  which is defined by
$$
\Delta _n(x)= \left|
\begin{array}{cccccc}
x&1&&&&\\
1&x&1&&&\\
&1&x&1&&\\
&&1&x& \ddots &\\
&&&\ddots &\ddots &1\\
&&&&1&x
\end{array} \right| _{n\times n}.$$
By Lemma \ref{lem1}, we have the following recursive relation:
$$\Delta _n(x)=x\Delta _{n-1}(x)-\Delta _{n-2}(x),
$$
with $\Delta _1(x)=x$ and $\Delta _2(x)=x^2-1$.  As such, it is concluded  from \cite{oxford} that
\begin{equation}\label{eq0}
\Delta _n(x)=U_n\left(\frac{x}{2}\right),
\end{equation}
where $U_n(x)$ is the $n$-th degree Chebyshev polynomial of second kind. Let $C_n(x)=|xI_n-A_n|$. A straightforward calculation shows that
 $C_1(x)=\sqrt{b}~\Delta_1((x-a)/\sqrt{b})$ and
 $C_2(x)=
 b~\Delta_2((x-a)/\sqrt{b})$.
 We claim that
 \begin{equation}\label{I}
 C_n(x)=
 (\sqrt{b})^n~\Delta_n\left(\frac{x-a}{\sqrt{b}}\right).
 \end{equation}
 We may assume that (\ref{I}) holds for any $k < n$. Lemma \ref{lem1} yields
\begin{align}
C_n(x)=&~(x-a) C_{n-1}(x)-bC_{n-2}(x) \nonumber\\
=&~(x-a)(\sqrt{b})^{n-1}
\Delta_{n-1}\left(\frac{x-a}{\sqrt{b}}\right)-(\sqrt{b})^{n}
\Delta_{n-2}\left(\frac{x-a}{\sqrt{b}}\right)\nonumber\\
=&~(\sqrt{b})^{n} \Delta_{n}\left(\frac{x-a}{\sqrt{b}}\right).\nonumber
\end{align}
Now, the proof follows from (\ref{eq0}).
\end{proof}
\end{lm}

\section{Homogeneous Monotonic P-polynomial Table Algebras}\label{p-tablealgebra}
Throughout this section, we find  the characters of homogeneous monotonic P-polynomial table algebras with finite dimension $d\geq 5$. To do so, we focus on calculating the eigenvalues of the first intersection matrix of  homogeneous monotonic P-polynomial table algebras which has a special form. In general, the first intersection matrix of P-polynomial table algebras is given by (\ref{matrix}), but for homogeneous monotonic P-polynomial  table algebras, the first
 intersection matrix  has a certain structure which is determined in the following theorem.

\begin{theorem}\label{theorem}
({\cite[Theorem $1.1$]{table}})  Let $(A,\mathbf{B})$ be a homogeneous monotonic P-polynomial table algebra with finite dimension $d$ and valency $k\geq 2$. Then the first intersection matrix of $(A,\mathbf{B})$ is a $(d+1) \times (d+1)$ matrix as follows.
$$
B_1=\left(
\begin{array}{cccccc}
0&1&&&&\\
k&\alpha &k-\alpha -1&&&\\
&k-\alpha -1&2(\alpha +1)-k&k-\alpha -1&&\\
&&\ddots &\ddots &\ddots &\\
&&&k-\alpha -1&2(\alpha +1)-k&k-\alpha -1\\
&&&&k-\alpha -1&\alpha +1
\end{array} \right),
$$
and one of the following items holds:
\begin{itemize}
\item[(i)]
$d\geq 5$ and $\alpha =(3k-6)/4$; or
\item[(ii)]
$d=4$ and  $(2k-4)/3\leq \alpha \leq (3k-6)/4$; or
\item[(iii)]
$d=3$ and  $(k-2)/2\leq \alpha \leq (3k-6)/4$; or
\item[(iv)]
$d=2$ and  $0\leq \alpha \leq k-2$.
\end{itemize}
\end{theorem}
We now find the eigenvalues of the first intersection matrix of homogeneous monotonic P-polynomial table algebras with finite dimension $d\geq 5$ in the following theorem.
Note that for finite dimension $d\leq 4$, the calculation of the eigenvalues of the first intersection matrix depends on the value of $\alpha$ is not too complicated.
\begin{theorem}\label{theorem2}
 Let $(A,\mathbf{B})$ be a homogeneous monotonic P-polynomial table algebra with $\mathbf{B}=\{x_0=1_A, x_1, \cdots , x_d\}$ and $d\geq 5$. Then  the eigenvalues of the first intersection matrix of $(A,\mathbf{B})$ are given by
$$\lambda _i=\frac{k+2}{2}\cos(\theta _i)+\frac{k-2}{2},~~~0\leq i\leq d,$$
where $k$ is the valency of $(A,\mathbf{B})$ and the $\theta _i$ are the roots of the following equation:
\begin{align}\label{eq7}
& (k+2)
\sin ((d+2)\theta)-4 \sin ((d+1)\theta)-2k\sin (d\theta)\nonumber\\
+&4\sin ((d-1)\theta)
+(k-2)\sin ((d-2)\theta)=0.\nonumber
\end{align}
\begin{proof}
The first intersection matrix of $(A,\mathbf{B})$ is given in Theorem \ref{theorem}. Letting $a=k-\alpha -1$, we can rewrite the  first intersection matrix as follows.
$$
B_1=
\left(
\begin{array}{cccccc}
2a-k&1&&&&\\
k&a-1 &a&&&\\
&a&0&a&&\\
&&\ddots &\ddots &\ddots &\\
&&&a&0&a\\
&&&&a&a
\end{array} \right)
+\left(2(\alpha +1)-k\right) ~\mathrm{I}_{_{(d+1) \times (d+1)}}.
$$
Denote the above tridiagonal matrix by $P$. Then the  eigenvalues of $B_1$ are equal to the eigenvalues of $P$ plus $(k-2)/2$. Let $\lambda$ be an eigenvalue of $P$ and $u=[u_1, \cdots, u_{d+1}]^t$  be the eigenvector corresponding to $\lambda$. Then we can consider the eigenvector $u$ as a sequence $\{u[i]\}_{i=0}^{\infty}$ with
$$
u[i]=\begin{cases}
u_i, ~~ \text{if} ~ i=1, \cdots, d+1 \\
0, ~~ \text{otherwise.}
\end{cases}
$$
Since $Pu=\lambda u$, we have
$$
\begin{cases}
u[0]=0\\
au[0]+au[2]=\lambda u[1]+(k-2a)u[1]+(a-1)u[2]\\
au[1]+au[3]=\lambda u[2]+(a-k)u[2]+(1-a)u[2]\\
au[2]+au[4]=\lambda u[3]\\
\vdots \\
au[d-1]+au[d+1]=\lambda u[d]\\
au[d]+au[d+2]=\lambda u[d+1]-au[d+1]\\
u[d+2]=0.
\end{cases}
$$
Consequently, we have the following equation:
\begin{equation}\label{equ1}
au[h+2]+au[h]=\lambda u[h+1]+f[h+1], ~~ h=0, 1, \cdots
\end{equation}
where
$$
f[h]=
\begin{cases}
(k-2a)u[1]+(a-1)u[2], ~~~\text{if}~ h=1\\
(a-k)u[1]+(1-a)u[2], ~~~~ \text{if}~ h=2\\
-au[d+1],~~~~~~~~~~~~~~~~~~~~~\text{if}~ h=d+1\\
0, ~~~~~~~~~~~~~~~~~~~~~~~~~~~~~
~~~~\text{otherwise}.
\end{cases}
$$
The $\textbf{z}$-transform of $f[h]$ is
\begin{align}\label{equ2}
F[z]=&\sum _{h=0}^{\infty}f[h]z^{-h}\nonumber\\
=&\left( (k-2a)u[1]+(a-1)u[2]\right)z^{-1}\nonumber\\
+&\left((a-k)u[1]+(1-a)u[2]\right)z^{-2}-au[d+1]z^{-(d+1)}.
\end{align}
From (\ref{pro3}), we calculate the $\textbf{z}$-transform of~(\ref{equ1}) which is
\begin{align}\label{eq3}
U(z)=&~\frac{1}{az^2-\lambda z+a}\left(F[z]+au[1]\right)z\nonumber\\
=&~\frac{1}{az^{-2}-\lambda z^{-1}+a}\left(F[z]+au[1]\right)z^{-1}.
\end{align}
Let $\displaystyle{X(z)= {1\over az^{-2}- \lambda z^{-1}+ a}}$. From (\ref{pro2}), we know that the inverse $\textbf{z}$-transform of $z^{-q}X(z)$
is $x[j-q]H[j-q]$, where
$$
H[x]=\begin{cases}
1,~~~~ \text{if}~ x\geq 0\\
0, ~~~~\text{if}~ x< 0.
\end{cases}
$$
From (\ref{equ2}) and (\ref{eq3}), we get
\begin{align}\label{eq4}
u[j]=&~au[1] x[j-1]H[j-1]+\left((k-2a)u[1]+(a-1)u[2]\right)x[j-2]H[j-2]\nonumber\\
+&\left((a-k)u[1]+(1-a)u[2]\right)x[j-3]H[j-3]\nonumber\\
-&au[n]x[j-d-2]H[j-d-2].
\end{align}
We shall now find $x[j]$ and plug it into (\ref{eq4}). Note that the eigenvalues of $B_1$ are real because $x_1=\overline{x_1}$, see  {\cite[Lemma $2.8$]{quotient}}, and therefore the coefficients of $az^{-2}-\lambda z^{-1}+a$ are all real and three cases may arise as follows.
\begin{itemize}
\item[(i)] If $az^{-2}-\lambda z^{-1}+a$ has two complex conjugate roots, then the roots are
$$\gamma _{\pm}=\frac{\lambda \pm \sqrt{\omega}}{2a},$$
where $\omega =\lambda ^2-4a^2<0$. We can write $\gamma _{\pm}=p\pm iq$,  where $p, q\in \mathbb{C}$ and $q\neq 0$. From $\gamma _+\gamma_ -=p^2+q^2=1$ and $\gamma _+ +\gamma _-=2p=\lambda /a$, it follows that
$$\gamma _{\pm}=\sqrt{p^2+q^2}\left(\cos (\theta)+i\sin (\theta)\right)=e^{\pm i\theta},~~~ \cos (\theta)=\frac{\lambda}{2a}$$
By  partial fractions decomposition of $X(z)$, we can write
$$X(z)=\frac{1}{\sqrt{\omega}}\left(\frac{1}{\gamma _ --z^{-1}}-\frac{1}{\gamma _+-z^{-1}}\right),$$
therefore,
\begin{equation}\label{eq5}
x[j]
=\frac{1}{\sqrt{\omega}}(\gamma _+^{j+1}-\gamma _-^{j+1})
=\frac{2i}{\sqrt{\omega}}\sin ((j+1)\theta).
\end{equation}
Note that $\{a^{-(n+1)}\}_{n=0}^{\infty}\leftrightarrow 1/(a-z^{-1})$. From (\ref{eq4}) and (\ref{eq5}), we have
\begin{align}\label{eq6}
u[j]=&\frac{2i}{\sqrt{\omega}}\Bigg[au[1]\sin(j\theta)H[j-1]\nonumber\\
+&\left((k-2a)u[1]+
(a-1)u[2]\right)\sin((j-1)\theta)H[j-2]\nonumber\\
+&\left((a-k)u[1]+(1-a)u[2]\right)\sin((j-2)\theta)H[j-3]\nonumber\\
-&au[n]\sin((j-d-1)\theta)H[j-d-2]\Bigg].
\end{align}
Setting  $j=d+2$ in (\ref{eq6}) yields
\begin{align}\label{eq16}
u[d+2]=&\frac{2i}{\sqrt{\omega}}\Bigg[au[1]\sin((d+2)\theta)\nonumber\\
+&\left((k-2a)u[1]+
(a-1)u[2]\right)\sin((d+1)\theta)\nonumber\\
+&\left((a-k)u[1]+(1-a)u[2]\right)\sin(d\theta)\nonumber\\
-&au[d+1]\sin(\theta)\Bigg].
\end{align}
Similarly, when $j= d+1$,
\begin{align}\label{eq17}
u[d+1]=&\frac{2i}{\sqrt{\omega}}\Bigg[au[1]\sin((d+1)\theta)\nonumber\\
+&\left((k-2a)u[1]+
(a-1)u[2]\right)\sin(d\theta)\nonumber\\
+&\left((a-k)u[1]+(1-a)u[2]\right)\sin((d-1)\theta)\Bigg].
\end{align}
Moreover, we know that
\begin{equation}\label{eq18}
\gamma _+ -\gamma _-=\frac{\sqrt{\omega}}{a}=2i\sin (\theta).
\end{equation}
From (\ref{eq16}), (\ref{eq17}) and (\ref{eq18}), we conclude that
\begin{align}
u[d+2]=&\frac{2i}{\sqrt{\omega}}\Bigg[au[1]\sin((d+2)\theta)\nonumber\\
+&\left((k-3a)u[1]+
(a-1)u[2]\right)\sin((d+1)\theta)\nonumber\\
+&\left((3a-2k)u[1]+(2-2a)u[2]\right)\sin(d\theta)\nonumber\\
-&\left((a-k)u[1]+(1-a)u[2]\right)\sin((d-1)\theta)\Bigg].
\end{align}
By considering $u[2]=(k-2a+\lambda)u[1]$, $a=(k+2)/4$ and $\lambda =2a\cos (\theta)$, we get
\begin{align}
u[d+2]=\frac{iu[1](k+2)}{8\sqrt{\omega}}~&\Bigg[ (k+2)
\sin ((d+2)\theta)\nonumber\\
-&4 \sin ((d+1)\theta)-2k\sin (d\theta)\nonumber\\
+&4\sin ((d-1)\theta)
+(k-2)\sin ((d-2)\theta)\Bigg].
\end{align}
Since $u[d+2]=0$, $\theta$ is the root of the following equation:
\begin{align}\label{eq7}
& (k+2)
\sin ((d+2)\theta)-4 \sin ((d+1)\theta)-2k\sin (d\theta)\nonumber\\
+&4\sin ((d-1)\theta)
+(k-2)\sin ((d-2)\theta)=0.
\end{align}
\item[(ii)]Suppose that $az^{-2}-\lambda z^{-1}+a$ has a repeated root. As such, $\omega =0$ and $\lambda =\pm 2a$. This means
\begin{align}\label{eq8}
x[j] &= \frac{1}{az^{-2}-\lambda z^{-1}+a}\nonumber\\
&=\frac{1}{a(z^{-2}-\frac{\lambda}{a}
z^{-1}+1)}\nonumber\\
 &=\frac{1}{a(z^{-1}\pm 1)^2}\nonumber\\
&=\frac{1}{a}(\pm 1)^{j+1}j.
\end{align}
Note that the last equality is obtained due to $\{na^{n}\}_{n=0}^{\infty}\leftrightarrow az/(z-a)^2$. From (\ref{eq4}) and (\ref{eq8}), we have
\begin{align}\label{eq9}
u[j]=&\frac{1}{|a|}\Bigg[au[1](\pm 1)^j(j-1)H[j-1]\nonumber\\
+&\left((k-2a)u[1]+(a-1)u[2]\right)(\pm 1)^{j-1}(j-2)
H[j-2]\nonumber\\
+&\left((a-k)u[1]+(1-a)u[2]\right)(\pm 1)^{j-2}(j-3)H[j-3]\nonumber\\
-&
au[d+1](\pm 1)^{j-d-1}(j-d-2)H[j-d-2]\Bigg].
\end{align}
If we set $j=d+2$ and $u[2]=(k-2a+\lambda)u[1]$ in (\ref{eq9}), we have
\begin{align}
&(d+1)\Bigg[a(\pm 1)^{d+2}+\left((k-2a)+(a-1)(k-2a+\lambda)\right)
(\pm 1)^{d+1}\nonumber\\
&+\left((a-k)+(1-a)(k-2a+\lambda)\right)(\pm1)^{d}\Bigg]\nonumber\\
&-\Bigg[\left((k-2a)+(a-1)(k-2a+\lambda)\right)
(\pm 1)^{d+1}\nonumber\\
&+2\left((a-k)+(1-a)(k-2a+\lambda)\right)(\pm1)^{d}\Bigg]=0.\nonumber
\end{align}
The above equation can be solved only for $k=2$ which means some eigenvalues of $P$ are $+2a=2$ or $-2a=-2$ for  $k=2$.
\item[(iii)] If $az^{-2}-\lambda z^{-1}+a$ has two real and distinct roots, then the roots are
$$\gamma _{\pm}=\frac{\lambda \pm \sqrt{\omega}}{2a},$$
where $\omega =\lambda ^2-4a^2>0$. Then
$$X(z)=\frac{1}{\sqrt{\omega}}\left(\frac{1}{\gamma _ --z^{-1}}-\frac{1}{\gamma _+-z^{-1}}\right),$$
and
\begin{align}
x[j]
=&~\frac{1}{\sqrt{\omega}}(\gamma _+^{j+1}-\gamma _-^{j+1}).
\end{align}
  From (\ref{eq4}), we find $u[j]$ as follows.
\begin{align}\label{eq10}
u[j]=&\frac{1}{\sqrt{\omega}}\Bigg[au[1](\gamma _+^j-\gamma _-^j)H[j-1]\nonumber\\
+&\left((k-2a)u[1]+(a-1)u[2]\right)(\gamma _+^{j-1}-\gamma _-^{j-1})H[j-2]\nonumber\\
+&\left((a-k)u[1]+(1-a)u[2]\right)(\gamma _+^{j-2}-\gamma _-^{j-2})H[j-3]\nonumber\\
-&au[d+1](\gamma _+^{j-d-1}-\gamma _-^{j-d-1})H[j-d-2]\Bigg].
\end{align}
Next, we set  $j=d+2$ and calculate $u[d+2]$:
\begin{align}
u[d+2]=&~\frac{1}{\sqrt{\omega}}\Bigg[au[1](\gamma _+
^{d+2}-\gamma _-^{d+2})\nonumber\\
+&\left((k-3a)u[1]+(a-1)u[2]\right)(\gamma _+^{d+1}-\gamma _-^{d+1})\nonumber\\
+&\left((3a-2k)u[1]+(2-2a)u[2]\right)(\gamma _+^{d}-\gamma _-^{d})\nonumber\\
-&\left((a-k)u[1]+(1-a)u[2]\right)(\gamma _+^{d-1}-\gamma _-^{d-1})\Bigg].\nonumber
\end{align}
Since $u[d+2]=0$, we have
\begin{align}
 & \Bigg[au[1]\gamma _+^3+\left((k-3a)u[1]+(a-1)u[2]\right)\gamma _+^2
\nonumber\\
&+\left((3a-2k)u[1]+(2-2a)u[2]\right)\gamma _+-\left((a-k)u[1]+(1-a)u[2]\right)\Bigg]\nonumber\\
&\Bigg[au[1]\gamma _-^3+
\left((k-3a)u[1]+(a-1)u[2]\right)\gamma _-^2
\nonumber\\
&+\left((3a-2k)u[1]+(2-2a)u[2]\right)\gamma _--\left((a-k)u[1]+(1-a)u[2]\right)\Bigg]^{-1}\nonumber\\
&=\left(\frac{\gamma _-}{\gamma _+} \right)^{d-1}.\nonumber
\end{align}
Given that $|\gamma _-/\gamma _+|<1$, the infinite series $\sum_{d} \gamma _-/\gamma _+$ is convergent. As such, if we consider the series $\sum_{d} b_d$, where $b_d$ is the left hand side of the above equality, then this series must be convergent which implies $\lim _{d\rightarrow \infty}b_d=0$. On the other hand, the term $b_d$ is not dependent
on $d$ and
is a constant, so $\sum_{d} b_d=0$. In other words, we have $a=0$,
which is a contradiction since $a=(k+2)/4$. This means item $(iii)$ does not occur and $az^{-2}-\lambda z^{-1}+a$  can not have two distinct real roots.
\end{itemize}

All in all, we conclude that the eigenvalues of $B_1$  can all be obtained from (\ref{eq7}), so the proof is completed.
\end{proof}
\end{theorem}
In the following theorem, we calculate the characters of homogeneous monotonic P-polynomial table algebras with finite dimension $d\geq 5$.
\begin{theorem}\label{theorem3}
Let $(A,\textbf{B})$  be  a homogeneous monotonic P-polynomial table algebra with  $\mathbf{B}=\{x_0=1_A, x_1, \cdots , x_d\}$ and $d\geq 5$. Then the characters of $(A,\textbf{B})$ are
$$p_0(j)=1, ~~~ 0 \leq j\leq d,$$
$$p_1(j)=\lambda _j, ~~~ 0 \leq j\leq d,$$
\begin{align}
p_i(j)=&\left( \frac{\sqrt{k+2}}{2}\right)^{i-4}\Bigg[\left(\lambda _j^2-\frac{3k-6}{4}\lambda_j-k\right)U _{i-2}\left(\frac{2\lambda _j-k+2}{2\sqrt{k+2}}\right)\nonumber\\
-&\left(\frac{\sqrt{k+2}}{2}\right)^3\lambda _jU _{i-3}\left(\frac{2\lambda _j-k+2}{2\sqrt{k+2}}\right)\Bigg],~~~ 0 \leq j\leq d,\nonumber
\end{align}
for $2\leq i \leq d$, where the $\lambda _j$ are the eigenvalues of the first intersection matrix of $(A,\mathbf{B})$ as given in Theorem \ref{theorem2} and $U_n$ is the $n$-th degree Chebyshev polynomial of second kind.

\begin{proof}
For each $i$, $0\leq i\leq d$, the $p _i(j)$, $0\leq j \leq d$, are equal to the eigenvalues of the $i$-th intersection matrix $B_i$.
Since $B_0=I_{d+1}$, we have $p _0(j) = 1$ for all $j$. Similarly, the $p _1(j)$ are equal to the eigenvalues of $B_1$ which are  calculated in Theorem \ref{theorem2}. To obtain the $p_i(j)$, $2\leq i\leq d$,
we must calculate the complex cofficient polynomial $\nu _i(x)$, where $x_i=\nu _i(x_1)$. Obviously, $\nu _1(x)=x$, and from ~(\ref{abc}) we get
$$x_{_1}x_{_1}=k+\alpha x_{_1}+(k-\alpha -1)x_{_2}\Rightarrow\nu _{_2}(x)=\frac{1}{k-\alpha -1}\left(x^2-\alpha x-k\right).$$
We claim that
\begin{equation}\label{nu}
\nu _{_i}(x)=\frac{1}{k-\alpha -1}\Bigg[\left(x-2(\alpha +1)+k\right)\nu _{_{i-1}}(x)-(k-\alpha -1)\nu _{_{i-2}}(x)\Bigg],~~~3\leq i\leq d.
\end{equation}
To prove this, we use induction on $i$. It is fairly straightforward using (\ref{abc}) to get
$$\nu _{_3}(x)=\frac{1}{k-\alpha -1}\Bigg[\left(x-2(\alpha +1)+k\right)\nu _{_2}(x)-(k-\alpha -1)\nu _1(x)\Bigg].$$
Now, we assume that (\ref{nu}) holds for $i<d$. From (\ref{abc}) and  the induction hypothesis, it is concluded that
$$\nu _{_d}(x)=\frac{1}{k-\alpha -1}\Bigg[\left(x-2(\alpha +1)+k\right)\nu _{_{d-1}}(x)-(k-\alpha -1)\nu _{_{d-2}}(x)\Bigg].$$

We now consider the following recursive relation:
\begin{equation}\label{recurs}
\varphi _n(x)=\left(x-2(\alpha +1)+k\right)\varphi _{n-1}(x)-(k-\alpha -1)\varphi _{n-2}(x);~~~~~~n>2,\nonumber
\end{equation}
with $\varphi _1(x)=(k-\alpha -1)x$ and $\varphi _2(x)=x^2-\alpha x-k$. Setting $\alpha =(3k-6)/4$ and Lemma \ref{lem1} yield
$$\varphi _n(x)=\left|
\begin{array}{cccccc}
\frac{k+2}{4}x&1&&&&\\
k&\frac{4}{k+2}x-\frac{3k-6}{k+2} &1&&&\\
&\frac{k+2}{4}&x-\frac{k-2}{2}&1&&\\
&&\ddots &\ddots &\ddots &\\
&&&\frac{k+2}{4}&x-\frac{k-2}{2}&1\\
&&&&\frac{k+2}{4}&x-\frac{k-2}{2}
\end{array} \right| _{n\times n}.$$
Laplace expansion gives
\begin{equation}\label{phi}
\varphi _n(x)=\left(x^2-\frac{3k-6}{4}x-k\right)D _{n-2}(x)-\left(\frac{k+2}{4}\right)^2xD _{n-3}(x),
\end{equation}
where
$D_n(x)$ is the characteristic polynomial of
$$
\left(
\begin{array}{ccccc}
\frac{k-2}{2}&1&&&\\
\frac{k+2}{4}&\frac{k-2}{2}&1&\\
&\frac{k+2}{4}&\frac{k-2}{2}&\ddots\\
&&\ddots &\ddots&1\\
&&&\frac{k+2}{4}&\frac{k-2}{2}
\end{array} \right)_{n\times n}.
$$
From Lemma \ref{lem2}, we have
\begin{equation}\label{D}
D_n(x)=\left( \frac{\sqrt{k+2}}{2}\right)^nU_n\left(\frac{2x-k+2}{2\sqrt{k+2}}\right).
\end{equation}
Finally, from (\ref{nu}), (\ref{phi}) and (\ref{D}) we conclude that
\begin{align}
\nu _{_i}(x)=~&\frac{1}{k-\alpha -1}\varphi _{_i}(x)\nonumber\\
=~&\left( \frac{\sqrt{k+2}}{2}\right)^{i-4}\Bigg[\left(x^2-\frac{3k-6}{4}x-k\right)U _{i-2}\left(\frac{2x-k+2}{2\sqrt{k+2}}\right)\nonumber\\
-~&\left(\frac{\sqrt{k+2}}{2}\right)^3xU _{i-3}\left(\frac{2x-k+2}{2\sqrt{k+2}}\right)\Bigg],\nonumber
\end{align}
for $2\leq i\leq d$. Due to (\ref{character}), the proof is now complete.
\end{proof}
\end{theorem}
\subsection{Two Classes of P-polynomial Association Schemes}
In general, the Bose-Mesner algebra of an association scheme is a table algebra. We can apply the characters of table algebras to calculate the eigenvalues of association schemes. For instance, we now give two classes of P-polynomial association schemes which are studied in \cite{perfect}. Note that the Bose-Mesner algebra of a P-polynomial association scheme is a monotonic P-polynomial table algebra, {\cite[Proposition  III.$1.2$]{Bannai}}.
\begin{example}\label{ex1}
Let  $(A,\textbf{B})$ be a homogeneous monotonic P-polynomial table algebra of valency $k=2$ and diameter $d\geq 5$.
Then from Theorem \ref{theorem}, the first intersection of $(A,\textbf{B})$ is
$$
B_1=\left(
\begin{array}{cccccc}
0&1&&&&\\
2&0 &1&&&\\
&1&0&1&&\\
&&\ddots &\ddots &\ddots &\\
&&&1&0&1\\
&&&&1&1
\end{array} \right)_{(d+1)\times (d+1)}.
$$
Now, we  calculate the characters of $(A,\textbf{B})$.  From Theorem \ref{theorem2},  we must find the roots of the following equation:
\begin{equation}\label{last}
    \sin ((d+2)\theta)- \sin ((d+1)\theta)-\sin (d\theta)
+\sin ((d-1)\theta)=0.
\end{equation}
It implies that
$$
\cos((d+1)\theta)\sin(\theta)
-\cos(d\theta)\sin(\theta)=0.
$$
We can assume that $\sin(\theta)\neq 0$ because otherwise $\theta =n\pi$, $(n\in \mathbb{N})$, and $\lambda=2a\cos(n\pi)=\pm 2a$ which is a contradiction. It is because (\ref{last}) is obtained from case $(i)$ in Theorem \ref{theorem2} in which we assume that $X(z)$ has two complex conjugate roots. This gives
$$p_1(j)=\lambda _j=2\cos\left(\frac{2j\pi}{2d+1}\right),~~~~0\leq j \leq d.$$
The other characters of $(A, \textbf{B})$ can also be found through Theorem~\ref{theorem3}. The $p_i(j)$, $(2\leq i\leq d)$, are
\begin{align}
p_i(j)=&\left(\lambda _j^2-2\right)U _{i-2}\left(\frac{\lambda _j}{2}\right)
-\lambda _j U _{i-3}\left(\frac{\lambda_j}{2}\right)\nonumber\\
=&~2\cos\left(\frac{ij\pi}{2d+1}\right),\nonumber
\end{align}
for $0\leq j \leq d$, where $U_i(x)$ and $T_i(x)$ are the $i$-th degree Chebyshev polynomials of second kind and first kind, respectively. The above equalities follow from the properties of Chebyshev polynomials  which can be found in
 \cite{oxford}.
\end{example}
In example \ref{ex1}, $(A,\textbf{B})$ is isomorphic to the Bose-Mesner algebra of the P-polynomial association scheme
$\mathcal{C}_1=(X,\{R_i\}_{0\leq i\leq d})$ of the ordinary $n$-gon where $n=2d+1$, see {\cite[Theorem  $1.5$]{perfect}}. In this case, the $n$-gon is a Moore graph with valency 2 and diameter $d$. Moore graphs are a class of distance-regular graphs which are introduced in {\cite[Section  $3.3$]{Bannai}}. Additionally, by calculating the eigenvalues of  $\mathcal{C}_1$, we are able to study the multiplicities and Krein parameters of  $\mathcal{C}_1$. As such, from {\cite[Theorem  III.$3.5$]{Bannai}},  the multiplicitie are  as follows.
$$m_0=1,~~~~~~ m_r=2,~~1\leq r \leq d.$$
Moreover, from {\cite[Theorem  III.$3.6$]{Bannai}}, the Krein parameters of $\mathcal{C}_1$ are obtained as shown in Table \ref{table2}. In Table \ref{table2}, each cell is a
column-vector and the $k$-component of the vector in row $i$
and column $j$ is $q_{ij}^k$.  As expected, since the Bose-Mesner algebra of $\mathcal{C}_1$ is a self-dual table algebra, the Krein parameters of $\mathcal{C}_1$ are $0$, 1, or 2 and  are the same as the intersection numbers.

\begin{example}\label{ex2}
Let $(A,\textbf{B})$ with $\mathbf{B}=\{x_0, x_1,\cdots ,x_d\}$ be a table algebra which is isomorphic to the Bose-Mesner algebra of the P-polynomial association scheme $\mathcal{C}_2$ of the ordinary $n$-gon, where $n=2d$. The first intersection matrix of $(A,\textbf{B})$ is
$$B_1= \left(
\begin{array}{cccccc}
0&1&&&&\\
2&0&1&&&\\
&1&0&\ddots &&\\
&&1&0 & 1 &\\
&&&\ddots &\ddots &2\\
&&&&1&0
\end{array} \right)_{(d+1) \times (d+1)},$$
see {\cite[Theorem  $1.5$]{perfect}}.
 $(A,\textbf{B})$  is not homogeneous because $f(x_i)=2$, $(0<i<d)$, and $f(x_d)=1$. See {\cite[Proposition  III.$1.2$]{Bannai}} for details on the calculation of $f(x_i)$. As a result, we can not directly apply Theorems \ref{theorem2} and \ref{theorem3} to calculate the characters of $(A,\textbf{B})$, but we apply some techniques for tridiagonal matrices and the argument which is used to prove Theorem \ref{theorem3} to obtain the characters. Let $H_{_{d+1}}(x)=|xI_{_{d+1}}-B_1|$ 
and $P_n(x)$ be a function as follows:
$$
P_n(x)= \left|
\begin{array}{cccccc}
x&2&&&&\\
1&x&1&&&\\
&1&x&1&&\\
&&1&x& \ddots &\\
&&&\ddots &\ddots &1\\
&&&&2&x
\end{array}\right|_{n\times n}.
$$
From Lemma \ref{lem1}, $H_{_{d+1}}(x)=P_{_{d+1}}(x)$ and hence from \cite{eigenvector}, the eigenvalues of $B_1$ are
$$p_1(j)=2\cos\left(\frac{j\pi}{d} \right), \,\,\  0\leq j \leq d. $$
In order to calculate the $p _i(j)$, $(2\leq i \leq d)$, we should obtain the complex cofficient polynomial $\nu_i(x)$, where $x_i=\nu_i(x_1)$. From the first intersection matrix of $(A,\mathbf{B})$, we can obtain
$$
\nu _i(x)=x\nu_{i-1}(x)-\nu_{i-2}(x), \,\ (\nu _1(x)=x, \nu_ 2(x)=x^2-2),
$$
for $3\leq i\leq d-1$ and
\begin{align}\label{nud}
&\nu _d(x)=\frac{1}{2}
\left(x\nu _{d-1}(x)-\nu_{d-2}(x)\right).
\end{align}
On the other hand, $\nu _i(x)$ is the same as $\varphi _i(x)$ in Theorem \ref{theorem3} for $k=2$, which implies that
\begin{equation}\label{nud-1}
\nu _i(x)=2T_i(\frac{x}{2}), ~~~1\leq i\leq d-1.
\end{equation}
Thus from (\ref{character}), we have
\begin{align}
p_i(j)=&
~2\cos\left(\frac{ij\pi}{d} \right),\nonumber
\end{align}
for $2\leq i\leq d-1$ and $ 0\leq j\leq d$.
Finally, we can calculate the $p_d(j)$ using (\ref{character}),(\ref{nud}) and (\ref{nud-1}) as follows.
\begin{align}
p_d(j)=&~
\cos(j\pi )=(-1)^j,\nonumber
\end{align}
for
$0\leq j \leq d$.
\end{example}
In example \ref{ex2}, the $n$-gon is a generalized Moore graph with valency 2, see {\cite[Section  $3.3$]{Bannai}} for more details. The multiplicities of the association scheme $\mathcal{C}_2$ in the above example are calculated as follows.
$$m_0=m_d=1,~~~~ m_r=2,~~ 1\leq r\leq d-1.$$
Moreover, the Krein parameters of $\mathcal{C}_2$
 are given in Table \ref{table1}.
\section{Concluding Remarks}\label{remarks}
In this paper, we use the $\textbf{z}$-transform concept along with techniques from linear algebra and matrix theory in order to calculate the characters of homogeneous monotonic P-polynomial table algebras with finite dimension $d\geq 5$.  Importantly, we calculate the eigenvalues of a special classes of tridiagonal matrices which may have applications in other fields. Next, we obtain the characters of homogeneous monotonic P-polynomial table algebras with finite dimension $d\geq 5$ in terms of Chebyshev polynomials. Finally, we apply our results to calculate the eigenvalues of two classes of P-polynomial association schemes which come from Moore graphs and generalized Moore graphs with valency 2.

\newpage

\begin{table}[ht]
\begin{center}
\caption{Krein Parameters of $\mathcal{C}_1$, ($r, s\in \{1, \cdots, d\}$).}
{\begin{tabular}{|c|c|c|}
\hline
\backslashbox{i}{j} & 0 & $\mathrm{s}$ \\
\hline
&&\\
& $q_{00}^0=1$ & $q_{0s}^0=0$ \\
0&&\\
&$q_{00}^w=0$&$
q_{0s}^w=\begin{cases}
1, & \text{if} ~~s=w\\
0 ,& \text{otherwise}
\end{cases}
$\\
&&\\
\hline
&&\\
&$q_{r0}^0=0$&$q_{rs}^0=\begin{cases}
2, & \text{if} ~~r=s\\
0,& \text{otherwise}
\end{cases}$ \\
r &&\\
& $q_{r0}^w=\begin{cases}
1, & \text{if} ~~r=w\\
0 ,& \text{otherwise}
\end{cases}$&$q_{rs}^w=0$\\
&&\\
\hline
\end{tabular}}
\label{table2}
\end{center}
\end{table}
\newpage
\begin{table}[h]
\begin{center}
\caption{Krein Parameters of $\mathcal{C}_2$, ($r, s, w\in \{1, \cdots, d-1\}$).}
{\begin{tabular}{|c|c|c|c|}
\hline
\backslashbox{i}{j} & 0 & $\mathrm{s}$ & $\mathrm{d}$ \\
\hline
&&&\\
& $q_{00}^0=1$ & $q_{0s}^0=0$ & ~~~~$q_{0d}^0=0$~~~~ \\
&&&\\
0&$q_{00}^w=0$&$
q_{0s}^w=\begin{cases}
1, & \text{if} ~~s=w\\
0 ,& \text{otherwise}
\end{cases}
$&$q_{0d}^w=0$\\
&&&\\
&$q_{00}^d=0$&$q_{0s}^d=0$&$q_{0d}^d=1$\\
&&&\\
\hline
&&&\\
&$q_{r0}^0=0$&$q_{rs}^0 =\begin{cases}
2, & \text{if} ~~r=s\\
0 ,& \text{otherwise}
\end{cases}$&$q_{rd}^0= 0$ \\
&&&\\
r&$q_{r0}^w=\begin{cases}
1, & \text{if} ~~r=w\\
0 ,& \text{otherwise}
\end{cases}$&$q_{rs}^w=0$&$q_{rd}^w=0$\\
&&&\\
&$q_{r0}^d=0$&$q_{rs}^d=0$&$q_{rd}^d=0$\\
&&&\\
\hline
&&&\\
&$q_{d0}^0=0$&$q_{ds}^0 =0$&$q_{dd}^0=1$ \\
&&&\\
d&$q_{d0}^w=0$&$q_{ds}^w=0$&$q_{dd}^w=0$\\
&&&\\
&$q_{d0}^d=1$&$q_{ds}^d=0$&$q_{dd}^d=0$\\
&&&\\
\hline
\end{tabular}}
\label{table1}
\end{center}
\end{table}
\vspace{2cm}


\begin{thebibliography}{99}

\bibitem{Bannai} E. Bannai, T. Ito, \textit{ Algebraic Combinatorics I: Association Schemes.} Menlo Park, CA: Benjamin/Cummings. 1984.

\bibitem{quotient} H. I. Blau, \textit{ Quotient structures in $C$-algebras}. J. Algebra $\mathbf{177}$ (1995), no. 1, 297--337.

\bibitem{table} H. I. Blau, R. J. Hein, {\em A Class of P-polynomial Table Algebras with and without Integer Multiplicities}. Comm. Algebra $\mathbf{42}$ (2014), 5387--5424.

\bibitem{Fib} N. D. Cahill, J. R. \'DEricco,  J. P. Spence, {\it Complex factorizations of the Fibonacci and Lucas numbers}, Fibonacci Quart. {\bf 41} (2003), no. 1, 13--19.

\bibitem{eigensystem} H. W. Chang, S. E. Liu, R. Burridge, {\it Exact eigensystems for some matrices arising from discretizations}, Linear Algebra Appl. {\bf 430} (2009), no. 4, 999--1006.


\bibitem{structure} G. Chen, B. Xu, \textit{Structures of commutative table algebras determined by their character tables and applicatons to finite groups.} Comm. Algebra \textbf{46}
(2018), no. 8, 3510--3519.


\bibitem{Gods} C. Godsile, {\em Association Schemes 1, Combinatorics, Optimization }. University of Waterloo
(2010).



\bibitem{oxford} L. Fox, I. B. Parker, {\it Chebyshev polynomials in numerical analysis}, Oxford University Press, London, 1986.

\bibitem{matlab} E. W. Kamen,  B. S. Heck, {\em Fundamentals of Signals and Systems: Using the Web and MATLAB}, 2007.

\bibitem{several} S. Kouachi, {\em Eigenvalues and eigenvectors of tridiagonal matrices}, Electron. J. Linear Algebra {\bf 15} (2006), 115--133.

\bibitem{system} J. Rimas, {\em Investigation of the dynamics of mutually synchronized systems}, Telecommun. Eng. {\bf 32} (1977), 68--79.


\bibitem{eigenvector} J. Rimas, {\em On computing of arbitrary positive integer powers for one type of tridiagonal matrices}, App. Math. Comput.  \textbf{161} (2005), no.
3, 1037--1040.


\bibitem{characters} B. Xu,  {\it Characters of table algebras and applications to association schemes}.  J. Combin. Theory Ser. A $\mathbf{115}$ (2008), no. 8, 1358--1373.

\bibitem{class} B. Xu, \textit{On a class of integral table algebras.} J. Algebra \textbf{178} (1995), no. 3, 760--781.

\bibitem{perfect} B. Xu, \textit{On perfect $P$-polynomial table algebras},
 Comm. Algebra \textbf{37}  (2009), no. 1, 120--153.

\bibitem{scheme} B. Xu, \textit{On P-polynomial table algebras and applications to association schemes.} Comm. Algebra \textbf{40} (2012), no. 6, 2171--2183.


\bibitem{Anal} A. J. Willms, {\em Analytic results for the eigenvalues of certain tridiagonal matrices}, SIAM J. Matrix Anal. Appl. $\mathbf{30}$ (2008), no. 2, 639--656.

\end{thebibliography}
\end{document}